\documentclass[11pt]{amsart}
\usepackage{amsmath,amscd,amssymb,amssymb}

\newtheorem{The}{Theorem}

\newtheorem{Def}{Definition}[section]
\newtheorem{Lem}[Def]{Lemma}

\newtheorem{Pro}[Def]{Proposition}
\newtheorem{Obs}[Def]{Observation}
\newtheorem{Ex}[Def]{Example}

\newtheorem{Cla}[Def]{Claim}

\newtheorem{Rem}[Def]{Remark}

\newenvironment{standard-model}[1]{\vspace{2ex} \noindent {\bf Standard models
                for conformal actions on {#1}:} \vspace{1ex}\\}{\\}

\newcommand {\la}{\lambda}
\newcommand {\ka}{\kappa}

\newcommand {\f}{\phi}
\newcommand {\F}{\Phi}

\newcommand {\id}{\mbox{id}}

\def \Fq{{\mathbb{F}}_q}

\newcommand {\Hf}{\mathcal{H}_{\f}}

\newcommand {\gb}{\mathbf{g}}
\newcommand {\cb}{\mathbf{c}}
\newcommand {\eb}{\mathbf{e}}

\newcommand {\ta}{\tau}
\newcommand {\tav}{\vec{\tau}}
\newcommand {\Tav}{\vec{\mathcal{T}}}

\newcommand {\vMM}{\vec{\mathcal{M}}}
\newcommand {\vTT}{\vec{\mathcal{T}}}
\newcommand {\vCC}{\vec{\mathcal{C}}}

\newcommand {\Gb}{\mathbf{G}}

\newcommand {\Cb}{\mathbf{C}}

\newcommand {\Tb}{\mathbf{T}}

\newcommand {\N}{\mathbb{N}}

\newcommand {\C}{\mathbb{Q}}
\newcommand {\Q}{\mathbb{Q}}
\newcommand {\HH}{\mathcal{H}}

\newcommand {\KK}{\mathcal{K}}
\newcommand {\NN}{\mathcal{N}}
\newcommand {\FF}{\mathcal{F}}
\newcommand {\CC}{\mathcal{C}}

\newcommand {\PP}{\mathcal{P}}

\newcommand {\OO}{\mathcal{O}}

\newcommand {\UU}{\mathcal{U}}

\newcommand {\TT}{\mathcal{T}}
\newcommand {\MM}{\mathcal{M}}

\newcommand {\half}{\frac{1}{2}}
\newcommand {\Span}{\mathrm{span}}
\newcommand {\Ker}{\mathrm{Ker}}

\newcommand {\rad}{\mathrm{rad}}

\newcommand {\coker}{\mathrm{coker}}

\newcommand {\Iso}{\mathrm{Iso}}
\newcommand {\Aut}{\mathrm{Aut}}
\newcommand {\Hom}{\mathrm{Hom}}
\newcommand {\End}{\mathrm{End}}

\newcommand {\Ob}{\mathrm{Ob}}

\newcommand {\GL}{\mathrm{GL}}

\newcommand {\Symp}{\mathrm{Symp}}
\newcommand {\Gr}{\mathrm{Gr}}

\newcommand {\Iff}{if and only if}

\font\tenolde=eufm10 at 10pt
 \font\sevenolde=eufm7
 \font\fiveolde=eufm5
 \newfam\oldefam
 \textfont\oldefam=\tenolde
 \scriptfont\oldefam=\sevenolde
 \scriptscriptfont\oldefam=\fiveolde

\newcommand{\embi}[3]{#1 \underset{#2}{\hookrightarrow} #3}

\newcommand{\ad}[2]{\Tb_{#2 \prec #1}}

\newcommand{\ade}[3]{\Tb_{#2 \underset{#3}{\hookrightarrow} #1}}

\newcommand{\au}[2]{\Tb_{#2 \succ #1}}

\newcommand{\aue}[3]{\Tb_{#2 \underset{#3}{\hookleftarrow} #1}}

\newcommand{\geomie}[4]{\Gb_{#3 \succ #2 \underset{#4}{\hookrightarrow} #1}}

\newcommand{\ceomie}[4]{\Cb_{#3 \succ #2 \underset{#4}{\hookrightarrow} #1}}

\newcommand{\cont}[3]{[#1 \!\prec \!\dot{#2}]_{#3}}

\newcommand{\avoid}[3]{[#1 \!\prec \!\dot{#2} \!\pitchfork\! #3]}

\title{Geometric representations of $\GL(n,R)$, cellular Hecke algebras and the embedding problem}
\author{Uri Bader and Uri Onn}

\thanks{Supported by Israel Science Foundation, ISF grant no.
100146.}

\begin{document}

\maketitle

\begin{abstract}
We study geometric representations of $GL(n,R)$ for a ring $R$. The
structure of the associated Hecke algebras is analyzed and shown to
be cellular. Multiplicities of the irreducible constituents of these
representations are linked to the embedding problem of pairs of
R-modules $x \subset y$.
\end{abstract}

\section{Introduction}

Let $R$ be a ring and let $F=R^n$ be the free module of rank $n$. In
this paper we study families of representations of
$G=GL(n,R)=\Aut_{R}(F)$ which arise from its action on the lattice
of submodules of $F$. More precisely, let $_R\MM$ denote the
category of finitely generated left $R$-modules, and let $\la \in
\Iso(_R\MM)$ be an isomorphism type of a submodule of $F$. Let
$X_{\la}=\Gr(\la,F)$ be the Grassmannian of submodules of type $\la$
in $F$. Let $\FF_{\la}$ be the vector space of $\Q$-valued functions
on $X_{\la}$. We define a family of representations of $G$
\begin{equation*}
\begin{split}  \rho_{\la}:G \longrightarrow \Aut_{\Q}(\FF_\la)
 \qquad \qquad \quad \\g \longmapsto [\rho_{\la}(g)f](x)=f(g^{-1}x). \end{split}
\end{equation*}
For each representation $\FF_{\la}$, let
$\HH_{\la}=\End_G(\FF_{\la})$ be the Hecke algebra associated to it.
The aim of this paper is two fold:

\medskip

First, we analyze the structure of the Hecke algebras
$\{\HH_{\la}\}$ for a distinguished family of $\la$'s and establish
their cellular structure, which is a consequence of the rich
underlying geometrical/combinatorial structure. This is a
generalization of \cite{hill} in which the ring in question is a
discrete valuation ring. Moreover, the technics used here are valid
in a much broader setting, when $F$ is an object in a category $\CC$
satisfying some axioms, and the group $G$ is $\Aut_{\CC}(F)$.

\medskip

Second, an equality is established between multiplicities of a
family of irreducible representations $\{\UU_{\mu}\}_{\mu \in
\Iso(_R\MM)}$ in the $\FF_\la$'s and the cardinality of
non-equivalent pairs $x \subset y$ of $R$-modules. As a result we
give a representation theoretic view on the embedding problem
\cite{RS,Schmidmeier}.

\subsection{Description of results}

To make things concrete, we explain our results in the special case
$R=\OO$, a discrete valuation ring with finite residue field. Let
$\OO$ be such a ring, let $\wp$ be its maximal ideal and let
$\OO_k=\OO/\wp^k$ for some $k \in \N$. By the principal divisors
theorem, any finite $\OO$-module is of the form $M=\oplus_{i=1}^{l}
\OO/\wp^{\la_i}$ with decreasing exponents and its isomorphism type
is the partition $\la=(\la_1,\ldots,\la_l)$. The set
$\TT=\Iso(_R\MM)$ of isomorphism types of $R$-modules carries a
natural structure of a poset, the partial order defined by: $\la \le
\nu$ if and only if a module of type $\la$ can be embedded in a
module of type $\nu$. Denote by $\{\la\hookrightarrow\nu\}$ the set
of arrow types with source of type $\la$ and range of type $\nu$,
and fix $\f=k^m=(k,\ldots,k)$ where $m \le n/2$. Then the following
theorem is a specialization to $R=\OO$ of Theorem
\ref{multiplicities}.
\begin{The}
There exists a collection of non-equivalent irreducible
$G$-representations $\bigl\{\UU_{\la}\bigr\}_{\la \le \f}$ such
that:
\begin{enumerate}
    \item $\FF_{\f}=\bigoplus_{\la \le \f}\UU_{\la}$.
    \item For every $\la,\nu \le \f$:
\[ \langle\UU_{\la},\FF_{\nu}\rangle=|\{\la \hookrightarrow \nu\}| \]
I.e., the multiplicity of $\UU_{\la}$ in $\FF_{\nu}$ is the number
of non-equivalent embeddings of a module of type $\la$ in a module
of type $\nu$. In particular $\UU_{\la}$ appears in $\FF_{\la}$ with
multiplicity one and does not appear in $\FF_{\mu}$ unless $\la \le
\nu$.
\end{enumerate}

\end{The}
For a general ring $R$, we isolate the properties of the type
$\f=k^m$ in the above theorem, and show that it is still valid for a
general ring for types satisfying these properties (see Definitions
\ref{sym} and \ref{symcoup}).

\medskip

Apart from the representations, another central object in this work
is the Hecke algebra $\HH_{\f}=\End_G(\FF_{\f})$. This algebra is
commutative and its irreducible representations correspond
bijectively to the irreducibles of $G$ which occur in $\FF_{\f}$. We
show that this algebra has a cellular structure, intimately related
to the poset of isomorphism types. One of the main results
concerning the Hecke algebra is the existence of two families of
geometrically defined ideals of $\HH_{\f}$, $\{\HH_{\f}^{\la}\}_{\la
\le \phi}$ and $\{\HH_{\f}^{\la^-}\}_{\la \le \phi}$, such that for
each $\la$, $\KK_{\la}=\HH_{\f}^{\la}/\HH_{\f}^{\la^-}$ is one
dimensional representation of $\HH_{\f}$, and

\begin{The}\qquad
\begin{enumerate}

\item $\{\KK_{\la}\}_{\la \le \phi}$ forms a complete set of
irreducible representations of $\HH_{\f}$.

\item As $\HH_{\f}$-modules,
\[
\HH_{\f}^{\la} \simeq \bigoplus_{\mu \le \la} \KK_{\mu}
\]

\end{enumerate}

\end{The}

\subsection{About the organization of this paper}

The paper is organized as follows. In Section \ref{gen} we give
notations, definitions and basic results on incidence algebras.
These are important ingredients in the dictionary between
combinatorial data, encoded in the poset structure, and the
algebraic or representation theoretic structure we seek to expose.
Section \ref{rep} is the core of this manuscript. There we define
and study the Hecke algebra and Hecke modules associated to the
various representations under consideration. In section
\ref{fourier} we open a discussion on the explicit computation of
the idempotents of the Hecke algebra. They are expressed in terms of
combinatorial invariants of the lattice of submodules. As an
illustration we give a new proof of the Fourier decomposition for
the case $R=\Fq$, the finite field with $q$ elements. In a sequel
paper \cite{BO2} we push these results further and give the Fourier
decomposition for the case $R=\OO$, a discrete valuation ring with
finite residue field. In Section \ref{general} we show how to
general most of the results in this paper to a much general
framework where the category $_R\MM$ is replaced with a general
category.

\subsection{Relevant references}
Some of the examples considered here (see \ref{examples}) have
already appeared in the literature. Using different methods,
Grassmannians of sets and vector spaces were studied in
\cite{delsarte-report}, \cite{dunkl} respectively. This study was
carried on in \cite{stanton} to somewhat more general Grassmannians
related to finite Chevalley groups. One of the main themes in these
papers concerns the relations between representations of the groups
and orthogonal polynomials. In a different perspective, the study of
Grassmannians of finite modules over discrete valuation rings was
initiated by G. Hill in \cite{hill}. Hill was motivated by the
classification of irreducible representations of $\GL_n(\OO)$, where
$\OO$ is the ring of integers of a non-archimedean local field. The
main theorem in \cite{hill} turns out to be a special case of
Theorem~\ref{cellalg}. Furthermore, the explicit computation of the
idempotents in this case will be carried out in a subsequent paper
\cite{BO2}, using a delicate version of the method developed in
Section \ref{fourier}.

\subsection{Acknowledgements}
We are most grateful to Amos Nevo for hosting and encouraging this
research. We thank Shai Haran for his interest, encouragement and
stimulating discussions and to Bernhard Keller for commenting on
categorical aspects. We also thank Claus Ringel and Markus
Schmidmeier for supplying us information on the embedding problem.

\section{The setup}\label{gen}
\subsection{Grassmannians}

Let $R$ be a ring and $\MM$ the category of finitely generated
$R$-modules. We denote by $\vMM$ the category whose objects are
injections in $\MM$, and the set of morphisms between two injections
$x \overset{i}{\hookrightarrow}  y$ and $x'
\overset{i'}{\hookrightarrow}  y'$ is the set of commutative squares
\begin{equation}\label{D1}
  \begin{matrix} x & \overset{i}{\hookrightarrow} &  y \\
                  \downarrow &     & \downarrow  \\
                 x' & \overset{i'}{\hookrightarrow}  & y'
  \end{matrix}
\end{equation}
There is a natural notion of subobjects in $\vMM$, namely, $i$ is a
subobject of $i'$, denoted $i \le i'$, if the vertical arrows in
\eqref{D1} are inclusions. Let $\TT=\pi_0(\MM)$ be the set of {\em
types} (isomorphism classes) in $\MM$; Similarly, let $\vTT$ be the
set of types in $\vMM$. $\TT$ and $\vTT$ have natural poset
structure: $\xi \le \eta$ if an object of type $\xi$ can be a
subobject of an object of type $\eta$. Let $\ta:\MM \rightarrow \TT$
and $\tav:\vMM \rightarrow \vTT$ be the {\em type maps}.

For $y \in \Ob(\MM)$ let $\MM_y$ be the lattice of submodules of
$y$. Denote $\eta=\tau(y)$. We define the \emph{Grassmannian of
submodules of type $\xi$ in $y$} and the {\em Grassmanniann of
submodules with embedding type $\iota$}
\begin{equation*}
\biggl({y \atop \xi}\biggr) = \{x \in \MM_y ~|~ \tau(x)=\xi\}~,
\qquad \qquad \biggl({y \atop \iota}\biggr) = \{x \in \MM_y ~|~
\tav(x \subset y)=\iota\}
\end{equation*}
their cardinality are denoted $\bigl({\eta \atop \xi}\bigr)$ and
 $\bigl({\eta \atop  \iota}\bigr)$
 respectively. There is a natural map $\Tav\rightarrow \TT^2$ taking a morphism
type to its source and range types. We denote by
$\{\xi\hookrightarrow\eta\}$ the preimage of $(\xi,\eta)$ under this
map.

\bigskip

Note that $\bigl({y \atop \xi}\bigr)$ and $\bigl({y \atop
\iota}\bigr)$ are an $\Aut(y)$-spaces. There are $|\{\xi
\hookrightarrow \eta\}|$ many $\Aut(y)$-orbits in $\bigl({y \atop
\xi}\bigr)$, thus
\[ \Bigl({\eta \atop \xi}\Bigr)=\sum_{\iota \in \{\xi
\hookrightarrow \eta\} }
   \Bigl({\eta \atop
   \iota}\Bigr).
\]
Of special interest is the case where the action is transitive. If
$y$ is a symmetric object, then the action of $\Aut(y)$ is
transitive.

\begin{Def}\label{sym}
We say that a type $\f$ is symmetric if for any $f$ and $f'$ of type
$\f$, $x \subseteq f$ and $x' \subseteq f'$ with an isomorphism $h:
x
 \overset{\sim}{\rightarrow} x'$, the following diagram could
be completed:

\begin{equation*}
  \begin{matrix} f & \overset{\sim}{\dashrightarrow} &  f' \\
                  \cup &     &  \cup \\
                 x & \overset{h}{\longrightarrow} & x'
  \end{matrix}
\end{equation*}
An object of symmetric type will be called a symmetric object.
\end{Def}

In the category of finite $\OO$-modules, a module $f$ is symmetric
if and only if it is free over $\OO/\text{Ann}(f)$. Furthermore, for
any pair $f \le F$ of such modules, of types $\f$ and $\F$
respectively, with $2\cdot\text{rank}(f) \le \text{rank}(F)$, the
couple $(\F,\f)$ is symmetric.

\begin{Def}\label{symcoup} Let $\f \le \F$ be symmetric types. The couple $(\F,\f)$ is
called a \emph{symmetric couple} if for every object $F$ of type
$\F$, every $ f \subseteq F$ of type $\f$ and for every $\la \le
\f$,

\begin{enumerate}
\item  There exist $f' \subseteq F$ such that $\ta(f')=\f$
and $\ta(f \wedge f')=\la$.

\item  For every $x \subseteq F$ such that $\tau(x) \le \f$
the square
\begin{equation}\tag{$\diamondsuit$} \label{cartesian}
  \begin{matrix} x & \subseteq & x \vee f \\
                  \cup &     &  \cup \\
                  x \wedge f& \subseteq &  f
  \end{matrix}
\end{equation}
is cartesian in $\MM$.

\end{enumerate}

\end{Def}

\begin{Rem} \qquad
\begin{enumerate}
\item The symbols $\wedge$ and $\vee$ stand for the meet and the join in
the lattice of submodules in $F$.
\item Explicitly, the second part of definition \eqref{symcoup} is
equivalent to $x \vee f = x+f \simeq x \oplus f / x \cap f$, with $x
\wedge f = x \cap f$ embedded in $x \oplus f$ via the diagonal
inclusion maps.
\end{enumerate}
 \end{Rem}

The properties of being symmetric and symmetric couple will be used
in the sequel to give a parametrization of couples of submodules up
to automorphisms of the ambient module, and therefore parameterize
bases of some Hecke algebras (see \S\ref{hecke}).

\subsection{Incidence Algebras}
Let $(\PP,\le)$ be a finite (more generally - locally finite)
poset and let $\PP^{(2)} \subseteq \PP \times \PP$ be the set of
all $(x,y)$ satisfying the relation $x \le y$. Following
\cite{foundation-i}, define the {\bf incidence algebra} of $\PP$
denoted by $I(\PP)$ to be the collection of $\Q$-valued functions
on $\PP^{(2)}$ with product given by:
\begin{equation*}
[f * g] (x,y) = \sum_{x \le z \le y} {f(x,z)g(z,y)}, \qquad x \le y
\end{equation*}
Endowed with this product, $I(\PP)$ becomes a unital associative
$\Q$-algebra with unit element given by Kronecker's delta function.
This algebra has two distinguished elements which will be used
frequently in the sequel: the {\bf zeta function} $\zeta_{\PP}= {\bf
1}_{\PP^{(2)}}$, and its inverse the {\bf M\"{o}bius function}
denoted by $\mu_{\PP}$ (cf. \cite{foundation-i}).

The M\"{o}bius function can be calculated inductively - for all $x
\in \PP$ set $\mu_{\PP}(x,x)=1$ and assuming $\mu_{\PP}(x,z)$ has
been calculated for all $z \in [x,y)$, set:

\begin{equation} \label{eq-mu}
\mu_{\PP}(x,y)=-\sum_{x \le z < y} \mu_{\PP}(x,z)
\end{equation}

  In the special case of $x=0$, we denote
  $\chi(y)=\mu_{\PP}(0,y)$. This notation is justified since
  equation~\ref{eq-mu} shows that $\mu_{\PP}(0,y)$ is the
  Euler characteristic of the flag complex associate to the poset
  $[0,y]$.

\begin{Obs} \label{mu}
  $\mu_{\PP}(x,y)$ depends merely on the poset $[x,y]$. More formally,
  \[ \mu_{\PP}(x,y)=\mu_{[x,y]}(0,1) \]
\end{Obs}

The vector space of $\Q$-valued functions on $\PP$, which we denote
by $V(\PP)$, has a natural structure of a module over the
incidence algebra $I(\PP)$. The action is given by:
\begin{equation*}
f \cdot v (x) = \sum_{y \ge x} f(x,y)v(y) \qquad f \in I(\PP), v
\in V(\PP)
\end{equation*}

Let $(\PP',\le')$ be another poset, and let $\tau:\PP\rightarrow
\PP'$ be a poset map. There exist a natural map from $V(\PP)$ to
$V(\PP')$ given by summation along fibers:
\begin{equation*}
 \ta_{\star}v(\xi)=\sum_{\ta(x)=\xi} v(x) \qquad v \in V(\PP), ~\xi \in
 \PP'
\end{equation*}
We would have liked to define a similar map from the incidence
algebras $I(\PP)$ to $I(\PP')$ which will be a homomorphism. In
general this is not possible. Nevertheless, such a morphism can be
defined on a subalgebra of $I(\PP)$. Denote
\[ J(\PP)=\big\{ f\in I(\PP)~|~\forall y_1,y_2
\in \PP \mbox{ with } \tau(y_1)=\tau(y_2), ~\forall x'\in \PP',~
\sum_{x \in \ta^{-1}(x')}f(x,y_1)= \sum_{x \in \ta^{-1}(x')}f(x,y_2)
\big\}
\]
It is easily verified that $J(\PP)$ is in fact a subalgebra of
$I(\PP)$. Define $\tau^{_{(1)}}_{\star}:J(\PP) \rightarrow I(\PP')$
by:
\begin{equation*}
[\tau^{_{(1)}}_{\star}f](x',y')=\sum_{x \in \ta^{-1}(x')}f(x,y)
\qquad x' \in \PP', ~ y \in
 \PP,~ \ta(y)=y'
\end{equation*}

\begin{Pro}\label{incidence} \quad
\begin{enumerate}
\item The map $\tau^{_{(1)}}_{\star}$ is a well defined
homomorphism of algebras.

\item The following diagram is commutative:
\begin{equation*}
  \begin{matrix} J(\PP) \times V(\PP) & \rightarrow &  V(\PP) \\
                  \tau^{_{(1)}}_{\star} \downarrow \qquad \ta_{\star}\downarrow \quad &
                  &  \ta_{\star}\downarrow \quad \\
                 I(\PP') \times V(\PP')  & \rightarrow  &  V(\PP')
  \end{matrix}
\end{equation*}
(the horizontal maps are the module maps)
\end{enumerate}

\end{Pro}

\begin{proof}
Direct calculation.
\end{proof}

Given $h\in J(\PP)$, we shall denote its image in $I(\PP')$ by
$\hat{h}=\tau^{_{(1)}}_{\star}(h)$.
 Let us introduce the following
property on the triple $(\PP,\PP',\tau)$.

\bigskip

\noindent ($\clubsuit$) ~~ {\em $\PP$ has a 0 element, and for every
$y_1,y_2\in\PP$ with $\tau(y_1)=\tau(y_2)$, the posets $[0,y_1]$ and
$[0,y_2]$ are isomorphic over $\PP'$, that is there is a poset
isomorphism $\pi:[0,y_1]\rightarrow [0,y_2]$ with
$\tau\circ\pi=\tau$.}

\bigskip

 Whenever property ($\clubsuit$) is satisfied, one easily sees that both
$\zeta_{\PP}$ and $\mu_{\PP}$ are in $J(\PP)$. For the two cases of
interest for us, the posets of submodules $\MM_z$ and the poset of
subinjections $\vMM_z$ with the type maps this is indeed the case.
We now specify the discussion to them.

\begin{Cla} \label{hat-mu}

  Fix an $R$-module $z$ and consider $\MM_z$, the poset of its submodules. For every $x \le y \in \MM_z$, the interval $[x,y]$ is
isomorphic as posets to the interval $[0,y/x]$. It follows that
$\mu_{\MM_z}(x,y)=\chi(y/x)$. Consider the type map $\tau: \MM_z
\rightarrow \TT$. Let $\alpha \le \beta$ be types, then

 \[ \hat{\mu}(\alpha,\beta)
     =\sum_{\iota \in \{ \alpha {\hookrightarrow} \beta \}}
     \sum_{~x: \vec{\tau}(x\subseteq y)=\iota}
     \mu(x,y)
     =\sum_{\iota \in \{ \alpha {\hookrightarrow} \beta \}~}
     \biggl({\beta
     \atop \iota }
     \biggr) \chi(\coker(\iota))
  \]
If furthermore $\beta$ is symmetric, then for every $x\le y$ of
  types $\alpha$ and $\beta$, the type of $y/x$ is the same, and
  will be denote $\beta/\alpha$.
  For a fixed module $y$ of type $\beta$ we have
  \[ \hat{\mu}(\alpha,\beta)=\sum_{x\in\tau^{-1}\alpha}\mu(x,y)
     =\sum_{x\in\tau^{-1}\alpha}\chi(\beta/\alpha)
     =\biggl({\beta \atop \alpha}\biggr)\chi(\beta/\alpha)
  \]

\end{Cla}

\begin{Ex} \label{mu-field}
  Assume $R$ is a finite field with $q$ elements. The flag complex associated to the poset
  $[0,x]$ is nothing but the Tits building associated to $x$, and
  it is known \cite{chen-rota,foundation-i} that
  \[ \chi(x)=(-1)^{\dim(x)}q^{\bigl({\dim(x) \atop 2}\bigr)}, \]
  therefore,
  \[ \mu_{\MM_z}(x,y)=(-1)^{\dim(y)-\dim(x)}q^{\bigl({\dim(y)-
     \dim(x) \atop 2}\bigr)},
  \]
and
  \[ \hat{\mu}(m,n)=(-1)^{n-m}q^{\bigl({n-m \atop
     2}\bigr)}\biggl({n \atop m}\biggr)_q.
  \]
\end{Ex}

\begin{Cla} Fix an $R$-module $z$ and consider $\vMM_z$, the poset of the
subinclusions of the identity map $z\rightarrow z$. Consider the
type map $\vec{\tau}:\vMM_z\rightarrow \vTT$. For a type $\alpha\le
\beta$ and an inclusion type
$\iota\in\{\alpha\hookrightarrow\beta\}$ we have
  \[ \hat{\mu}(\iota,\beta)
     =\biggl({\beta \atop \iota}
     \biggr)\chi(\coker(\iota))
  \]
\end{Cla}

\section{Representations}\label{rep}
\subsection{Representations arising from Grassmannians}

Let $F$ be a symmetric module in $_R\MM$ of type $\F$ (recall
Definition \ref{sym}). For each type $\la \in \TT$ let
$X_{\la}=\bigl({F \atop \la}\big)$ be the Grassmannian of submodules
of $F$ of type $\la$. Let $G=\Aut_{R}(F)$. As we mentioned in
\S\ref{gen}, $F$ being symmetric implies that $X_{\la}$ is a
homogeneous $G$-space. Let $\FF_{\la}$ stand for $\Q$-valued
functions on $X_{\la}$. Then the $\FF_{\la}$'s become a family of
representations of $G$:
\begin{equation*}
\begin{split}  \rho_{\la}:G \longrightarrow \Aut_{\Q}(\FF_\la)
 \qquad \qquad \quad \\g \longmapsto [\rho_{\la}(g)f](x)=f(g^{-1}x) \end{split}
\end{equation*}

\noindent $\FF_{\la}$ is equipped with the standard $G$-invariant
inner product:
\begin{equation*}
(f,g)=\sum_{x \in X_{\la}}f(x)g(x) \qquad \forall f,g \in
\FF_{\la}
\end{equation*}
Fix once and for all $\f \le \F$ which forms a symmetric couple
(Definition \ref{symcoup}). Recall that $\{\la\hookrightarrow\mu\}$
stands for the set of arrow types with source of type $\la$ and
range of type $\mu$. The main theorem in this section is:

\begin{The}\label{multiplicities} Let $F$ be a module in $_R\MM$ of symmetric type $\F$, let $\f$ be a type such that $\f \le \F$ is a
symmetric couple, and let $G=\Aut_R(F)$. Then, there exists a
collection of non-equivalent irreducible $G$-representations
$\bigl\{\UU_{\la}\bigr\}_{\la \le \f}$ such that:
\begin{enumerate}
    \item $\FF_{\f}=\bigoplus_{\la \le \f}\UU_{\la}$.
    \item For every $\la,\nu \le \f$:
\[ \langle\UU_{\la},\FF_{\nu}\rangle=|\{\la \hookrightarrow \nu\}| \]
I.e., the multiplicity of $\UU_{\la}$ in $\FF_{\nu}$ is the number
of non-equivalent embeddings of a module of type $\la$ in a module
of type $\nu$. In particular $\UU_{\la}$ appears in $\FF_{\la}$ with
multiplicity one and does not appear in $\FF_{\nu}$ unless $\la \le
\nu$.
\end{enumerate}

\end{The}

\begin{proof}
  Postponed to \S~\ref{proofs}.
\end{proof}

\subsection{An equivalence of categories}\label{G-Hecke}

Let $G$ be a group and denote by $\MM_G$ the category of its finite
dimensional representations. For any $V \in \MM_G$ we can slice
$\MM_G$ and look at the full subcategory $\MM_{G,V}$ which consists
of representations generated by irreducibles appearing in $V$. Let
$\HH_V=\End_G(V)$ be the \emph{Hecke algebra} associated with $V$
and let $\MM_{\HH_V}$ be the category of finitely generated left
$\HH_V$-modules.

\begin{Pro} \label{equivalence}
$\MM_{G,V}$ and $\MM_{\HH_V}$ are equivalent categories. The
irreducible representations in $V$ are in one to one
correspondence with the irreducible $\HH_V$-modules.
\end{Pro}

\begin{proof}
The functors:
\begin{align*}
   &\MM_{G,V} \longrightarrow  \MM_{\HH_V} &  \qquad  & \MM_{\HH_V}
   \longrightarrow\MM_{G,V} \\
  & \qquad U \longmapsto \Hom_G(V,U) & \qquad & \quad ~M \longmapsto V \otimes_{\HH_V}
  M \\
\end{align*}
forms an equivalence of categories.
\end{proof}

For $\la,\mu \in \TT$ set
$\NN_{\mu,\la}=\Hom_G(\FF_{\la},\FF_{\mu})$ and
$\HH_\la=\NN_{\la,\la}$. Composition gives paring:
\[ \NN_{\la,\mu} \times \NN_{\mu,\nu} \longrightarrow
\NN_{\la,\nu}\]

\noindent which turns $\HH_{\la}$ into an algebra and
$\NN_{\mu,\la}$ into an $\HH_{\la}$-$\HH_{\mu}$-bimodule. Our main
goal is to study $\HH_{\f}$ and its modules. It turns out that it is
fruitful to have a slightly broader picture.

\subsection{Notations for intertwining operators}

\begin{itemize}
\item Given an embedding type $\la \underset{i}{\hookrightarrow}
\mu$, attach to it the following operators (they are transpose of
each other):
 \begin{align*} \aue{\la}{\mu}{i}: &\FF_{\la} \longrightarrow \FF_{\mu}
\\  &~~h \longmapsto \aue{\la}{\mu}{i}h(y)=\sum_{ x: \tav(x \subseteq y)=i} h(x),  \qquad y \in X_{\mu} \\
 \ade{\mu}{\la}{i}: &\FF_{\mu} \longrightarrow \FF_{\la}
\\  &~~h \longmapsto \ade{\mu}{\la}{i}h (x)=\sum_{  y: \tav(x \subseteq y)=i } h(y),  \qquad x \in X_{\la}
\end{align*}

\item Given $\la \le \mu$ we can define the maps:
\begin{align*}
\au{\la}{\mu}=\sum_{\la \underset{i}{\hookrightarrow} \mu}
\aue{\la}{\mu}{i} && \ad{\mu}{\la}=\sum_{\la
\underset{i}{\hookrightarrow} \mu} \ade{\mu}{\la}{i}
\end{align*}
i.e. averaging over all submodules (supmodules) of same type
regardless of the embedding. If $\mu=\f$ is a symmetric type these
sums consist of a unique summand each, as there is unique embedding.
Then the two notations degenerate: $\au{\la}{\f}=\Tb_{\f
\hookleftarrow \la}$ and $\ad{\f}{\la}=\Tb_{\la \hookrightarrow
\f}$.

\end{itemize}

\noindent In order to minimize confusion, we follow the rule that
whenever an operator is labeled with a diagram, it acts from the
space indexed by the \emph{right} type of the diagram to the space
indexed by the \emph{left} type.

\subsection{The Hecke algebra $\Hf$}\label{hecke}

\subsubsection{The geometric basis of $\Hf$} \label{geoH}

Let $X_{\f} \times_G X_{\f}$ be the equivalence classes of the
diagonal action of $G$ on $X_{\f} \times X_{\f}$.
The set $X_{\f}
\times_G X_{\f}$ has a natural parametrization:
\[ X_{\f} \times_G X_{\f} \simeq \{ \la | \la \le \f\} = [0,\f]
\]
This isomorphism is given by the ($G$-invariant map) $X_{\f} \times
X_{\f}\rightarrow [0,\f]$, given by $(x,y)\mapsto \tau(x\wedge y)$.
This map is onto by definition \ref{symcoup} part (1) and is
one-to-one by part (2). Indeed, for $(x,y)$ and $(x',y')$ in $X_{\f}
\times X_{\f}$, if $\tau(x\wedge y)=\tau(x'\wedge y')$ then the
following diagram of inclusions could be completed to a commutative
one using vertical isomorphisms:
\begin{equation}\tag{D1}\label{diag}
\begin{matrix} x & \supseteq &  x \wedge y  & \subseteq & y \\
                  |\wr &     &   |\wr  &          & |\wr    \\
                  x' & \supseteq &  x' \wedge y'  & \subseteq & y' \end{matrix}
\end{equation}
By the cartesianity of the diagram \eqref{cartesian} in Definition
\ref{symcoup}, these isomorphisms extend to $x \vee y \simeq x' \vee
y'$. Finally, since $F$ is symmetric this isomorphism extends to an
automorphism of $F$. By viewing functions as integration kernels
$\Hf$ is seen to be isomorphic to $\FF(X_{\f}\times_G X_{\f})$,
which in turn is isomorphic to $\FF([0,\f])$ (by the above
discussion). Composing these isomorphisms we obtain the vector space
isomorphism
\begin{align*}
\FF([0,\f]) &\overset{\sim}{\rightarrow}\Hf \\
                          \mathbf{1}_{\la} &\mapsto \gb_{\la}
\end{align*}
where
\begin{equation*} [\gb_{\la}h](x)=\sum_{\{y \in X_{\f}|
\tau(x \wedge y)=\la\}}h(y) \qquad x \in X_{\f}, ~h \in \FF_{\f}
\end{equation*}
\begin{Rem}
 It is useful to view $\gb$ as a function
$\la\mapsto\gb_{\la}$, that is as an element of
$V([0,\f])\otimes\Hf$, that is the incidence module over the
subposet $[0,\f]$ in $\TT$ with coefficients in $\Hf$.
$V([0,\f])\otimes\Hf$ is a module over $I([0,\f])\otimes\Hf$ by
extension of scalars from $\Q$ to $\Hf$.
\end{Rem}

\begin{Pro}
  $\Hf$ is a semisimple commutative $\Q$-algebra of dimension
  $|[0,\f]|$.
\end{Pro}
\begin{proof}
The only thing we have to show is commutativity. This follows from
Gelfand's trick: transposition, which is an anti-isomorphism, is
realized by the flip map on $X_{\f}\times_G X_{\f}$, which is a
trivial operation.
\end{proof}

\subsubsection{The cellular structure of $\Hf$}

Define $\cb_{\la}=\au{\la}{\f} \ad{\f}{\la}$. Observe that up to
constant $\cb_{\la}$ is simply a composition of averaging
operators
    $\FF_{\f} \rightarrow \FF_{\la} \rightarrow \FF_{\f}$:

    \[ [\cb_{\la}h](x)=\sum_{X_{\la} \ni y \subseteq x} \sum_{X_{\f} \ni z \supset y} {h(z)}
    \qquad \qquad (h\in \FF_{\f}, ~x \in X_{\f})\]
Collecting terms according to intersection types and using the
definition of $\gb_{\la}$ gives:

\begin{equation}\tag{\textbf{c-g}} \label{c-g}
\cb_{\la}=\sum_{\ka \ge \la} \biggl({\ka \atop
\la}\biggr)\gb_{\ka}
\end{equation}

where $\bigl({\ka \atop \la}\bigr)=\hat{\zeta}(\la,\ka)$ is the
number of submodules of type $\la$ inside a module of type $\ka$.
Viewing $\cb$ as an element of $V([0,\f])\otimes \Hf$,
equation~\eqref{c-g} becomes $\cb=\hat{\zeta}\gb$. Therefore
$\gb=\hat{\mu}\cb$, and it follows that $\{\cb_{\la}\}_{\la\le\f}$
forms a basis of $\Hf$. We call this basis the {\em cellular basis}
of $\Hf$. For each $\la \le \f$ set:
\[\Hf^{\la}=\Span_{\C}\{\cb_{\la'} | \la' \le \la \}
\qquad \qquad \Hf^{\la^{-}}=\Span_{\C}\{\cb_{\la'} | \la' < \la \}
\]

\begin{The}[Cellular ideal structure]\label{cellalg} For every $\la,\mu \le
\f$:

\begin{enumerate}
\item $\Hf^{\la}$ and $\Hf^{\la^-}$ are ideals. In particular, they
contain a unit when considered as rings.

\item $\HH_{\f}^{\la} \cdot \HH_{\f}^{\mu} = \HH_{\f}^{\la} \cap
\HH_{\f}^{\mu}$. If $\TT$ is a lattice they
 both equal to $\HH_{\f}^{\la \wedge \mu}$.

\item $\bigl\{\KK_{\la}=\Hf^{\la}/\Hf^{\la^-}\bigr\}_{\la \le \f}$
is a complete set of inequivalent irreducible $\Hf$-modules.

\item As a $\HH_{\f}$-modules,
\[ \HH_{\f}^{\la}\simeq \bigoplus_{\mu\leq\la} \KK_{\mu} \]

\end{enumerate}

\end{The}

\begin{proof}
Postponed to \S\ref{proofs}.
\end{proof}

Note that Theorem~\ref{cellalg} can be used to characterize the
irreducible representations of $\Hf$: $\KK_{\la}$ is the
\emph{unique} representation which is \textbf{annihilated} by all
$\bigl\{\Hf^{\mu}\bigr\}_{\mu < \la}$ and \textbf{not annihilated}
by $\Hf^{\la}$. In view of the dictionary between representations
of $\Hf$ and representations of $G$ (section \ref{G-Hecke}), we
can now label the representations of $G$ which occur in $\FF_{\f}$
by: $\UU_{\la} \leftrightarrow \KK_{\la}$.

Moreover, by the definition of $\cb_{\la}$ as the composition
$\au{\la}{\f} \ad{\f}{\la}$, the annihilation criterion above
translates to the fact that $\UU_{\la}$ occurs in $\FF_{\la}$ and
does not occur in $\FF_{\mu}$ for $\la \nleq \mu$. The exact
multiplicities in which the $\UU_{\la}$'s appear in $\FF_{\mu}$ for
arbitrary $\mu$ is the subject of the next section.

\subsection{The Hecke modules $\NN_{\f,\nu}$}

In this section we consider the various bases, and the cellular
structure of the $\Hf$-module $\NN_{\f,\nu}$, for $\nu\le\f$. To
begin with, we give a simple Lemma.

\begin{Lem}\label{symcoup-gen}
Let $\la \le \nu \le \f$ and $\la \overset{i}{\hookrightarrow}
\nu$ an embedding type. Then there exist $x,f \subseteq F$, such
that $\tau(f)=\f$ and $\tav(x \wedge f \subseteq x)=i$.
\end{Lem}

\begin{proof}
 The proof is an exercise with the definitions. As $\f$ is symmetric, there exist $f,f'$ of type $\f$ with
  $\tau(f\wedge f')=\la$. Denote $y=f\wedge f'$.
  Fix a submodule $x'\le f'$ of type $\nu$.
  There exist also $y''\le x''$ with
  $\vec{\tau}(y''\le x'')=i$. Since $F$ is symmetric there is $g\in
  G$ with $gx''=x'\le f'$. Denote $gy''$ by $y'$. Since $f'$ is
  symmetric, there is some $h\in\Aut(f')$ such that $hy'=y$. We
  are done by letting $x=hx'$.
\end{proof}

\subsubsection{The geometric basis of $\NN_{\f,\nu}$}

We analyze the structure of $X_{\f} \times_G X_{\nu}$ following a
similar line to the analysis of $X_{\f} \times_G X_{\f}$, given in
\S\ref{geoH}.

The set $X_{\f} \times_G X_{\nu}$ has a natural parametrization:
\[ X_{\f} \times_G X_{\nu} \simeq \{\la\overset{i}{\hookrightarrow}
   \nu|\la\le\nu\}=[0,\nu \overset{\text{id}}{\rightarrow} \nu] \subset \vTT
\]
This isomorphism is given by the ($G$-invariant map) $X_{\f} \times
X_{\nu}\rightarrow [0,\nu \overset{\text{id}}{\rightarrow} \nu]$,
given by $(x,y)\mapsto \tav(x\wedge y \subseteq x)$. This map is an
isomorphism by the discussion in \S\ref{geoH} together with
Lemma~\ref{symcoup-gen}. $\NN_{\f,\nu}$ is isomorphic to
$\FF(X_{\f}\times_G X_{\nu})$, which in turn is isomorphic to
$\FF([0,\nu \overset{\text{id}}{\rightarrow} \nu])$. Composing these
isomorphisms we obtain the vector space isomorphism
\begin{align*}
V([0,\nu \overset{\text{id}}{\rightarrow} \nu])
&\overset{\sim}{\rightarrow}\NN_{\f,\nu} \\
 \mathbf{1}_{i} &\mapsto \Gb_{i}=\geomie{\nu}{\la}{\f}{i}
\end{align*}
where
\begin{equation*} [\Gb_i h](x)=\sum_{\{y \in X_{\nu}|
\tav(x \wedge y \subseteq x )=i\}}h(y) \qquad x \in X_{\f}, ~h \in
\FF_{\nu}
\end{equation*}
\begin{Rem} It is useful to view $\Gb$ as a function
$i\mapsto\Gb_{i}$, that is as an element of $V([0,\nu
\overset{\text{id}}{\rightarrow} \nu])\otimes\NN_{\f,\nu}$.
\end{Rem}

\subsubsection{The cellular structure of $\NN_{\f,\nu}$}

For an embedding type $\la \overset{i}{\hookrightarrow} \nu$, define
$\Cb_{i}=\ceomie{\nu}{\la}{\f}{i}= \au{\la}{\f} \ade{\nu}{\la}{i}$.
In analogy with (\textbf{c-g}), we have:

\begin{equation}\tag{\textbf{C-G}}\label{C-G}
\Cb_i=\sum_{(\embi{\eta'}{i'}{\la}) \ge (\embi{\eta}{i}{\la})}
\biggl({\embi{\eta'}{i'}{\la} \atop
\embi{\eta}{i}{\la}}\biggr)\Gb_{i'}
\end{equation}
where $\Bigl({\embi{\eta'}{i'}{\la} \atop
\embi{\eta}{i}{\la}}\Bigr)=\hat{\zeta}(i,i')$. Viewing $\Cb$ as an
element of $V([0,\nu \overset{\text{id}}{\rightarrow} \nu])\otimes
\Hf$, equation~\eqref{C-G} becomes $\Cb=\hat{\zeta}\Gb$. Therefore
$\Gb=\hat{\mu}\Cb$, and it follows that
$\{\Cb_{i}\}_{\embi{\la}{i}{\nu}}$ forms a basis of $\NN_{\f,\nu}$.
We call this basis the {\em cellular basis} of $\NN_{\f,\nu}$. For
each $\la \le \nu \le \f$ set:

\[\NN^{\la}_{\f,\nu}=\Span_{\C}\{ \Gb_i ~|~ \la'
\le \la, ~\la' \overset{i}{\hookrightarrow} \nu\} \]

\[ \NN^{\la^-}_{\f,\nu}=\Span_{\C}\{
\Gb_i | \la' < \la, ~\la' \overset{i}{\hookrightarrow} \nu\}
\]

\begin{The}[Cellular submodule structure]\label{cellmod} For every $\la \le \nu
\le\f$:

\begin{enumerate}
\item $\NN^{\la}_{\f,\nu}$ and $\NN^{\la^-}_{\f,\nu}$ are
submodules of $\NN_{\f,\nu}$.

\item
For every $\mu\le\f$,
$\HH_{\f}^{\mu} \NN_{\f,\nu}^{\la} =
    \NN_{\f,\nu}^{\mu} \cap \NN_{\f,\nu}^{\la}
$.
In particular,
$\Hf^{\la} \NN_{\f,\nu} = \NN_{\f,\nu}^{\la}$ and
$\Hf^{\la^-}\NN_{\f,\nu}=\NN^{\la^-}_{\f,\nu}$.

\item $\NN^{\la}_{\nu,\f}/\NN^{\la^-}_{\nu,\f}$ is isomorphic to
the $\KK_{\la}$-isotypic submodule of $\NN_{\f,\nu}$.

\item As a $\HH_{\f}$-modules,
\[ \NN_{\f,\nu}^{\la}\simeq \bigoplus_{\mu\leq\la} (\KK_{\mu})^{
   |\{\mu\hookrightarrow \nu\}|}
\]

\end{enumerate}

\end{The}
\begin{proof}
Postponed to \S\ref{proofs}.
\end{proof}

\subsection{Proofs}\label{proofs}

\subsubsection{Some Lemmas}

Denote by $\cont{\mu}{\la}{\f}$  the number of submodules of type
$\la$ which contain a given submodule of type $\mu$ and are
contained in a given object of type $\f$ (when $\f$ is symmetric
this is a well defined quantity).

\begin{Lem}\label{compo} Let $\mu \le \la \le \eta,\f$ be types,
assume $\f$ is symmetric. Let $i:\mu \hookrightarrow \la$ and
$j:\la \hookrightarrow \eta$ be given types of embeddings. Then
\begin{align*}
\mathbf{1.}~ &~\aue{\la}{\eta}{j}\aue{\mu}{\la}{i} \in
\Span_{\N}\{\aue{\mu}{\eta}{k}|k \} \quad &\mathbf{1'}.~
&~\ade{\la}{\mu}{i} \ade{\eta}{\la}{j} \in
\Span_{\N}\{\ade{\eta}{\mu}{k}| k \} \\ \mathbf{2.}~
&~\au{\la}{\f}\aue{\mu}{\la}{i} \in \N\cdot
 \au{\mu}{\f} & \mathbf{2'}.~ &~\ade{\la}{\mu}{i}\ad{\f}{\la} \in \N\cdot
 \ad{\f}{\mu} \\
\mathbf{3.}~ &~\au{\la}{\f}\au{\mu}{\la}= \cont{\mu}{\la}{\f}
    \au{\mu}{\f} & \mathbf{3'}.~
    &~\ad{\la}{\mu}\ad{\f}{\la}=\cont{\mu}{\la}{\f}\ad{\f}{\mu}
\end{align*}

\end{Lem}

\begin{proof} \textbf{1.} We fix $x_0 \in X_{\mu}$ and denote its characteristic function by
$\delta_{x_0} \in \FF_{\mu}$. Since this is a cyclic vector for
the representation, everything is determined by the action on this
element:
\begin{equation}\tag{$\star$}\label{aaa}
\begin{split}
 \aue{\la}{\eta}{j}\aue{\mu}{\la}{i} \delta_{x_0}(z) &= \sum_{\{y|y \underset{j}{\hookrightarrow}
z\}} \sum_{\{x|x \underset{i}{\hookrightarrow} y\}}
\delta_{x_0}(x) = \sum_{\{y|y \underset{j}{\hookrightarrow} z\}}
{\mathbf{1}}_{\{y'|x_0 \hookrightarrow y'\}}(y) \\ &= \sum_{\mu
\underset{k}{\hookrightarrow}
\eta} \#\biggl\{y\biggr|\begin{smallmatrix}  x_0 &\subseteq &y& \subseteq &z \\
\wr& & \wr && \wr \\ \mu &\underset{i}{\hookrightarrow}& \la
&\underset{j}{\hookrightarrow} & \eta
 \end{smallmatrix} \biggr\} \aue{\mu}{\eta}{k} \delta_{x_0}(z)
\end{split}
\end{equation}
In particular $\aue{\la}{\eta}{j}\aue{\mu}{\la}{i} \in
{\Span}_{\N}\{\aue{\mu}{\eta}{k}~|~k \}$.

\noindent \textbf{2.} Set $\eta=\f$ in the first part. Then
\eqref{aaa} reduces to:
\begin{equation} \tag{$\star\star$}\label{bbb}
\au{\la}{\f}\aue{\mu}{\la}{i}= \#\biggl\{y\biggr|\begin{smallmatrix}  x_0 &\subseteq &y& \subseteq &z \\
\wr& & \wr && \wr \\ \mu &\underset{i}{\hookrightarrow}& \la
&\hookrightarrow & \f
 \end{smallmatrix} \biggr\} \au{\mu}{\f}
\end{equation}

\noindent \textbf{3.} Follows from equation~(\ref{bbb}), as
\[ \sum_{i} \#\biggl\{y\biggr|\begin{smallmatrix}  x_0 &\subseteq &y& \subseteq &z \\
\wr& & \wr && \wr \\ \mu &\underset{i}{\hookrightarrow}& \la
&\hookrightarrow & \f
 \end{smallmatrix} \biggr\}=\cont{\mu}{\la}{\f}
\]

\noindent The other three assertions follow by transposition.
\end{proof}

\begin{Lem} \label{onto}
 Let $\mu\leq\la,\omega\le\f$ be types, and assume that the couple
 $(\F,\omega)$ is symmetric. Then the map
 \[ \NN_{\omega,\la}\rightarrow \HH_{\f},~~~C\mapsto \au{\omega}{\f}\circ C
    \circ\ad{\f}{\la}
 \]
 maps $\NN_{\omega,\la}^{\mu}$ onto $\HH_{\f}^{\mu}$, and it maps
 $\NN_{\f,\la}^{\mu^-}$ onto $\HH_{\f}^{\mu^-}$.
\end{Lem}

\begin{proof}
Follows immediately from Lemma~\ref{compo}, part (2), using
cellular bases.
\end{proof}

\begin{Lem}\label{ideals}
$\NN^{\la}_{\f,\nu}$ and $\NN^{\la^-}_{\f,\nu}$ are
$\Hf$-submodules of $\NN_{\f,\nu}$ $\forall \la \in [0,\nu]$.
\end{Lem}
\begin{proof}
Let $\mu \le \f$, $\la' \le \la$ and $i: \la' \hookrightarrow
\nu$.
\begin{align*}
\cb_{\mu} \cdot \Cb_i & = (\au{\mu}{\f} \ad{\f}{\mu}
\au{\la'}{\f})
\ade{\nu}{\la'}{i}  \\
& \in \Span_{\C}\{ \Cb_j | \eta \le \la',~\embi{\eta}{j}{\la'}
\} \circ \ade{\nu}{\la'}{i} &&  \\
&=\Span_{\C}\{\au{\eta}{\f}
\ade{\la'}{\eta}{j} \ade{\nu}{\la'}{i} | \eta \le \la',~\embi{\eta}{j}{\la'}\} &&\\
&\subseteq \Span_{\C}\{\au{\eta}{\f} \ade{\nu}{\eta}{k} | \eta \le
\nu,~ \embi{\eta}{k}{\nu}\} \\
&= \Span_{\C}\{\Cb_k | \eta \le \nu,~ \embi{\eta}{k}{\nu} \} =
\NN_{\f,\nu}^{\la} &&
\end{align*}
(the inclusion follows from Lemma \ref{compo})

\noindent This implies that $\NN_{\nu,\f}^{\la}$ is a module since
the $\cb_{\mu}$'s generate $\HH_{\f}$ and the $\Cb_i$'s generate
$\NN^{\la}_{\f,\nu}$. Finally $\NN^{\la^-}_{\f,\nu}$ is a
submodule as well, since $\NN^{\la^-}_{\f,\nu}=\sum_{\la' <
\la}\NN^{\la'}_{\f,\nu}$.
\end{proof}


\subsubsection{Proof of Theorem~\ref{cellalg}}

\begin{proof}
By Lemma \ref{ideals}, for every $\la\le\f$, $\Hf^{\la}$ is an
ideal (substitute $\nu=\f$). It follows that for every $\la\le\f$,
$\Hf^{\la^-}$ is an ideal as well, by: $\Hf^{\la^-}=\sum_{\la' <
\la}\Hf^{\la'}$. This proves part (1). To prove part (2), observe
first that $\Hf^{\la} \cdot \Hf^{\mu} \subseteq \HH_{\f}^{\la}
\cap \HH_{\f}^{\mu}$. The opposite inclusion follows from the fact
that each ideal (considered as a ring) contains a unit, namely the
sum of its minimal idempotents. Therefore:
\begin{equation*}
\HH_{\f}^{\la} \cdot \HH_{\f}^{\mu} \supseteq (\HH_{\f}^{\la} \cap
\HH_{\f}^{\mu}) \cdot (\HH_{\f}^{\la} \cap \HH_{\f}^{\mu}) =
\HH_{\f}^{\la} \cap \HH_{\f}^{\mu}
\end{equation*}
Part (3) follows from the fact that $\{\KK_{\la}\}_{\la \le \f}$
is a collection of $\dim(\Hf)$ distinct one dimensional
$\Hf$-representations (they are indeed distinct by the
characterization of $\KK_{\la}$ given after
Theorem~\ref{cellalg}). Part (4) is proved by induction with
respect to the partial order on $\TT$.
\end{proof}

\subsubsection{Proof of Theorem~\ref{cellmod}} Before proving the
theorem we shall need one more lemma.
\begin{Lem}\label{diamond}
Let $\la\le\nu\le\f$ and $\mu\le\f$ be types.
\begin{enumerate}
\item If $\mu \geq \la$ then
  \[ \HH_{\f}^{\mu} \NN_{\f,\nu}^{\la} =
     \NN_{\f,\nu}^{\la}
  \]
\item If $\mu \ngeq \la$ then
\[ \HH_{\f}^{\mu} \NN_{\f,\nu}^{\la} \subseteq
   \NN_{\f,\nu}^{\la^-}
\]
\end{enumerate}
\end{Lem}

\begin{proof}

\noindent {\bf Part 1:}
By Lemma~\ref{ideals}, $\HH_{\f}^{\mu}
\NN_{\f,\nu}^{\la} \subseteq \NN_{\f,\nu}^{\la}$. We argue to show
the reverse inclusion. Assume (by induction with respect to the
partial order on $\TT$) that for every $\la'<\la$,
$\HH_{\f}^{\mu}\NN^{\la'}_{\f,\nu}=\NN_{\f,\nu}^{\la'}$. Then
$\HH_{\f}^{\mu}\NN^{\la}_{\f,\nu}\supseteq\NN_{\f,\nu}^{\la^-}$,
and we are left to show that for every inclusion type $\la
\underset{i}{\hookrightarrow} \nu$, $\ceomie{\nu}{\la}{\f}{i} \in
\HH_{\f}^{\mu}\NN^{\la}_{\f,\nu}$. We first consider the case
$\nu=\la$. Assume
$\HH_{\f}^{\mu}\NN^{\la}_{\f,\la}=\NN_{\f,\la}^{\la^-}$. Composing
on the right with $\ad{\f}{\la}$ we get, using Lemma~\ref{onto},
that $\HH_{\f}^{\mu}\HH_{\f}^{\la}=\HH_{\f}^{\la^-}$, which is an
absurd, as $\HH_{\f}^{\la}$ has a unit. Therefore
$\HH_{\f}^{\mu}\NN^{\la}_{\f,\la}=\NN_{\f,\la}^{\la}$. In
particular we get that
$\au{\la}{\f}\in\HH_{\f}^{\mu}\NN^{\la}_{\f,\la}$. The general
case follows by Lemma~\ref{compo}, because
\[ \ceomie{\nu}{\la}{\f}{i} = \au{\la}{\f}\circ\ade{\nu}{\la}{i}
   \in \HH_{\f}^{\mu}\NN^{\la}_{\f,\la}\ade{\nu}{\la}{i} \subseteq
   \HH_{\f}^{\mu}\NN_{\f,\la}^{\nu}
\]
and we proved $\HH_{\f}^{\mu}\NN^{\la}_{\f,\nu}=
\NN_{\f,\nu}^{\la}$.

\noindent {\bf Part 2:} By Lemma~\ref{ideals}, $\Hf^{\mu}
\NN_{\f,\nu}^{\la^-} \subseteq \NN_{\f,\nu}^{\la^-}$. We have
therefore reduced the proof to showing that $\cb_{\mu'} \cdot
\ceomie{\nu}{\la}{\f}{i} = \au{\mu'}{\f} \ad{\f}{\mu'}
\au{\la}{\f} \ade{\nu}{\la}{i} \in \NN_{\f,\nu}^{\la^-}$ for all
$\mu' \le \mu$ and $i: \la \hookrightarrow \nu$. Express the map
$\au{\mu'}{\f} \ad{\f}{\mu'} \au{\la}{\f}$ in the cellular basis
of $\NN_{\f,\la}$:
\begin{equation} \label{jkl}
\au{\mu'}{\f} \ad{\f}{\mu'} \au{\la}{\f} = \sum_{\la'
\underset{j}{\hookrightarrow} \la} a_{\la'
\underset{j}{\hookrightarrow} \la} \au{\la'}{\f}
\ade{\la}{\la'}{j}
\end{equation}
\noindent We claim that $a_{\la \underset{id}{\hookrightarrow}
\la}=0$. Indeed, assume this was not the situation. Compose the
map $\ad{\f}{\la}$ on the right of both sides of
equation~(\ref{jkl}), and present it with respect to the cellular
basis of $\HH_{\f}$. Using Lemma~\ref{compo}, it is clear from the
presentation of the r.h.s. that the coefficient of $\cb_{\la}$ is
non-trivial. On the other hand, the l.h.s. becomes $\cb_{\mu'}
\cdot \cb_{\la}$. The latter is in $\HH_{\f}^{\la^-}$ by the
assumption $ \mu \ngeq \la$, using Theorem~\ref{cellalg}. This is
a contradiction.

Composing with $\ade{\nu}{\la}{i}$ on the right of both sides of
equation~(\ref{jkl}) gives:
\begin{equation*}\begin{split}
\au{\mu'}{\f} \ad{\f}{\mu'} \au{\la}{\f}\ade{\nu}{\la}{i} &=
\sum_{\la' \underset{j}{\hookrightarrow} \la, \la' < \la} a_{\la'
\underset{j}{\hookrightarrow} \la} \au{\la'}{\f}
(\ade{\la}{\la'}{j} \ade{\nu}{\la}{i})
\\ &= \sum_{\la' \underset{k}{\hookrightarrow}
\nu, \la' < \la} b_{\la' \underset{k}{\hookrightarrow} \nu}
\au{\la'}{\f} \ade{\nu}{\la'}{k} \in \NN_{\f,\nu}^{\la^-}
\end{split}\end{equation*}
\end{proof}

\begin{proof}[proof of Theorem~\ref{cellmod}]


The main part of the theorem is part (2).
Indeed, part (1) immidately follows from part (2).
Also, assuming part (2),
it follows that for every $\mu\le \f$, the $\KK_{\la}$-isotypic
component of $\NN_{\f,\nu}$ is
\[ \KK_{\la}\otimes_{\HH_{\f}} \NN_{\f,\nu} \simeq
   (\HH_{\f}^{\la}/ \HH_{\f}^{\la^-}) \otimes_{\HH_{\f}} \NN_{\f,\nu}
   \simeq \HH_{\f}^{\la} \NN_{\f,\nu}/ \HH_{\f}^{\la^-}
   \NN_{\f,\nu} \simeq \NN_{\f,\nu}^{\la}/ \NN_{\f,\nu}^{\la^-}
\]
and part (3) follows as well.
Finally, using cellular
bases, we see that $\dim (\NN_{\f,\nu}^{\la}/
\NN_{\f,\nu}^{\la^-})=|\{\la\hookrightarrow \nu \}|$, and part (4)
follows.

\smallskip

We proceed to the proof of part (2). By part (1)
lemma~\ref{diamond}, for every $\alpha\le\la,\mu$,
\[ \HH_{\f}^{\mu} \NN_{\f,\nu}^{\la} \supset
\HH_{\f,\nu}^{\alpha} \NN_{\f,\nu}^{\alpha}
= \NN_{\f,\nu}^{\alpha} \]
hence
\[ \HH_{\f}^{\mu} \NN_{\f,\nu}^{\la} \supset
\sum_{\alpha\le\la,\mu} \NN_{\f,\nu}^{\alpha}
= \NN_{\f,\nu}^{\mu} \cap \NN_{\f,\nu}^{\la}
 \]
Also we have that $\HH_{\f}^{\mu} \NN_{\f,\nu}^{\la} \subset
\HH_{\f} \NN_{\f,\nu}^{\la} = \NN_{\f,\nu}^{\la} $. We are left to
show that $\HH_{\f}^{\mu} \NN_{\f,\nu}^{\la} \subset
\NN_{\f,\nu}^{\mu}$. This follows by induction (with respect to
$\la$) on the poset $\TT$. Indeed, assume that for every $\la'<
\la$, $\HH_{\f}^{\mu} \NN_{\f,\nu}^{\la'} \subset
\NN_{\f,\nu}^{\mu}$. If $\mu \ge \la$ then we are done by part (1)
of lemma~\ref{diamond}. If $\mu \ngeq \la$, then recalling that
 $\HH_{\f}^{\mu}$ has a unit, and using part (2) of
lemma~\ref{diamond},
$$ \HH_{\f}^{\mu} \NN_{\f,\nu}^{\la} =
\HH_{\f}^{\mu}\left(\HH_{\f}^{\mu} \NN_{\f,\nu}^{\la}\right)
\subset
\HH_{\f}^{\mu}\left(\NN_{\f,\nu}^{\la-}\right)
=
\HH_{\f}^{\mu}\left(\sum_{\la'<\la}\NN_{\f,\nu}^{\la'}\right)
\subset
\NN_{\f,\nu}^{\mu}.$$

\end{proof}

\subsubsection{Proof of Theorem \ref{multiplicities}}

\begin{proof}
  By Proposition~\ref{equivalence}, The modules $\KK_{\la}$
  correspond to the irreducible representations of $G$,
  $\UU_{\la}=\FF_{\f}\otimes\KK_{\la}$. Part (2) of
  Theorem~\ref{multiplicities} now follows from
  Theorem~\ref{cellmod}, part (4). Part (1) is a special case of
  part (2), as $\f$ is symmetric.
\end{proof}


\section{Towards a Fourier decomposition} \label{fourier}

Throughout this section fix a symmetric type $\F$, and a module $F$
of type $\F$.

\subsection{A counting principle}\label{inc-exc}

Let $X$ be a set and assume we are given a map $\varphi:X
\rightarrow \MM_F$. Define the following elements in $V(\MM_F)$:
\begin{align*}
&s_{\varphi}(y) = |\{x \in X | \varphi(x)=y\}| \\
&t_{\varphi}(y) = |\{x \in X | \varphi(x) \supseteq y\}|
\end{align*}
Clearly $t_{\varphi}=\zeta \cdot s_{\varphi}$. By multiplying both
sides by $\mu$, using proposition \ref{incidence} applied to the map
$\tau:\MM_F \rightarrow \TT$, we get $\hat{s}_{\varphi}=\hat{\mu}
\cdot \hat{t}_{\varphi}$. Observing that $0 \in \MM_F$ is the
\emph{unique} element above $0 \in \TT$, we obtain an
inclusion-exclusion type formula which will be useful in the sequel:
\begin{equation} \label{basic}
s_{\varphi}(0)=\hat{s}_{\varphi}(0)=\hat{\mu} \cdot
\hat{t}_{\varphi}(0)=\sum_{\alpha}\chi(\alpha)\hat{t}_{\varphi}(\alpha)
\end{equation}

Recall (Lemma~\ref{compo}) that $[\alpha \prec \dot{\ka}]$ is
the number of submodule of type $\ka$ which contain a given submodule of type
$\alpha$. Let $x_{\ka}$ and $x_{\omega}$ be disjoint submodules of
$F$ of types $\ka$ and $\omega$ (i.e $x_{\ka}\wedge
x_{\omega}=0)$. We introduce another convenient notation:
$\avoid{\omega}{\beta}{\ka}_{\F}$ (or simply $\avoid{\omega}{\beta}{\ka}$)
will denote the number of
submodules of $F$ of type $\beta$ which contain $x_{\omega}$, and
are disjoint from $x_{\ka}$.

 Define \(
X=\{ x\subset F~|~\tau(x)=\beta,~x_{\omega}\subset x\} \), and
$\varphi: X\rightarrow \MM_F$ by $\varphi(x)=x\wedge x_{\ka}$.
Then, clearly,
\[ s_{\varphi}(0)=\avoid{\omega}{\beta}{\ka} \]
It is also easy to see that\[ \hat{t}_{\varphi}(\alpha)=\cont{\omega \oplus
   \alpha}{\beta}{}\Big({\ka \atop \alpha}\Big)
\]
Indeed,
\[ \hat{t}_{\varphi}(\alpha)=
\sum_{\{y\subset F~|~\tau(y)=\alpha\}} t_{\varphi}(y) =
\left|\bigcup_{\{y\subset x_{\ka}~|~\tau(y)=\alpha\}} \{x \subset F~|~\tau(x)=\beta,~y,x_{\omega}\subset x\}\right|
=\Big({\ka \atop \alpha}\Big)
\cont{\omega \oplus \alpha}{\beta}{}
\]
Thus, by equation~\eqref{basic},
\begin{equation} \label{specific}
\avoid{\omega}{\beta}{\ka}_{}=\sum_{\alpha}
\chi(\alpha)\cont{\omega \oplus \alpha}{\beta}{}\Big({\ka \atop
\alpha}\Big)
\end{equation}
(Observe that the right hand side is independent of the choice of
$x_{\omega}$ and $x_{\ka}$, hence so is the left hand side, and
the notation $\avoid{\omega}{\beta}{\ka}$ is justified).

 \begin{Ex}[The case of a field] \label{final}
  When $R$ is a field and $U,V$ are vector spaces, one has an embedding
  $\Hom(V,U)\hookrightarrow V\oplus U$, given by the {\bf graph}
  of a transformation. Its image consists of those subspaces
  intersecting $U$ trivially. If the field has $q$ elements,
  \[ \avoid{m}{l}{k}_n=\avoid{0}{(l-m)}{k}_{n-m}=
  \biggl( { n-m-k \atop l-m-k} \biggr)_{\! q}\cdot q^{k(l-m-k)}
  \]
  Thus, we get, using example~\ref{mu-field},
  \[
  \biggl( { n-m-k \atop l-m-k} \biggr)_{\! q}\cdot q^{k(l-m-k)} =
  \sum_i (-1)^i q^{\left( { i \atop 2} \right)}
  \biggl( { n-m-i \atop l-m-i} \biggr)_{\! q}
  \biggl( { k \atop i} \biggr)_{\! q}
  \]
  (compare with \cite{chen-rota})
 \end{Ex}

\subsection{Computing some matrix coefficients}

We wish to compute the idempotents of $\HH_{\f}$ explicitly. Since
the cellular structure must agree with the idempotent
decomposition, we already know that there exist a lower triangular
matrix $A_{\la\ka}$ such that

\begin{equation}\tag{\textbf{c-e}}
\cb_{\la}=\sum_{\ka \le \la} A_{\la \ka}\eb_{\ka}
\end{equation}
where $\eb_{\la}$ is the idempotent in $\HH_{\f}$ corresponding to
its irreducible representation $\KK_{\la}$. We have already seen
that the transition matrix from the geometric basis to the cellular
basis depends only on geometric invariants of the lattice of
submodules (relation (\textbf{c-g}) above) in a very simple way. In
some situations, we are able to give similar interpretation to
various matrix coefficients of $A$. This is the main theme of this
section.

Let $\ka\le\omega\le\F$ be types and assume $(\F,\omega)$ is a
symmetric couple. Assume that $\ka$ satisfies the following duality
axiom:

For a module $x$ of type $\ka$, and every type $\alpha\le \ka$,
\begin{equation} \tag{\textbf{duality}} \label{duality}
|\{y~|~y\le x,~\tau(y)=\alpha\}|=|\{y~|~y\le
x,~\tau(x/y)=\alpha\}|
\end{equation}

\noindent {\bf Remark:} By the principal divisor theorem, it is
easy to see that every finite module over a principal ideal domain
satisfies the duality axiom.

Recall that $[\ka \prec \dot{\omega}]_{\f}$ is the number of
submodules of type $\omega$ which contain a given submodule of
type $\ka$ and are contained in a given module of type $\f$ (that
is $[\ka \prec \dot{\omega}]_{\f}=[\ka \prec \dot{\omega}
\pitchfork 0 ]_{\f}$).

\begin{The}\label{cell-idem}
Under the above assumptions, $A_{\omega\ka} =
\cont{\ka}{\omega}{\f} \avoid{\omega}{\f}{\ka}_{\F}$.
\end{The}

Before proving Theorem~\ref{cell-idem} we state a simple Lemma.

\begin{Lem}
Let $\theta\underset{i}{\hookrightarrow}\ka$ be a map type. Let $w$
and $x$ be modules of types $\omega$ and $\ka$, such that
$\vec{\tau}(x\wedge w\le x)=i$. The following diagram is cartesian.
\begin{equation*}
  \begin{array}{rcl} w & \hookrightarrow &  w\oplus\coker(i) \\
                  \uparrow &     &  \uparrow  \\
                 w\wedge x & \underset{i}{\hookrightarrow} & x
  \end{array}
\end{equation*}
where
\begin{itemize}
\item $\coker(i)$ is the type of $x/(w\wedge x)$.
\item The map $w\hookrightarrow w\oplus\coker(i)$ is given by
$\id\oplus 0$.
\item The map $x\hookrightarrow
w\oplus\coker(i)$ is given $a\oplus b$ where $a$ is an embedding,
and $b$ is the natural projection.
\end{itemize}
\end{Lem}

\begin{proof}
  By the fact $\omega$ is symmetric, the embedding $w\wedge
  x\hookrightarrow w$ can be taken to be $-a\circ i$.
  One easily sees that
  \[ \Ker\bigl(\id\oplus b: w\oplus x\rightarrow  w\oplus\coker(i)\bigr)
     = \bigl((-a\circ i)\oplus i\bigr)(w\wedge x)
  \]
\end{proof}

\begin{proof}[Proof of Theorem~\ref{cell-idem}:]
Our strategy is to analyze the multiplication in the algebra with
respect to the cellular basis. Let $B_{\omega\mu}^{\nu}$ be
multiplication table with respect to the cellular basis:
\begin{equation}\label{mtable}
 \cb_{\omega} \cdot \cb_{\ka} = \sum_{\nu \le \omega \wedge \ka}
 B_{\omega\ka}^{\nu} \cb_{\nu}
 \end{equation}
Observe that $B^{\ka}_{\omega\ka}=A_{\omega\ka}$ for $\ka \le
\omega$.

\noindent Substituting $\cb_{\eta}= \au{\eta}{\f}\ad{\f}{\eta}$ in
equation \eqref{mtable} and using parts {\bf 3} and {\bf 3'} of
Lemma~\ref{compo} give (assume $\ka \le \omega$):

\begin{equation}\label{ddd}\begin{split}
 \au{\omega}{\f}\Bigl(&\ad{\f}{\omega}\au{\ka}{\f}\Bigr)\ad{\f}{\ka} \\
 &=\au{\omega}{\f} \Bigl(\sum_{\nu \le \ka}
\frac{B_{\omega\ka}^{\nu}}{\cont{\nu}{\omega}{\f}
\cont{\nu}{\ka}{\f}}
\au{\nu}{\omega}\ad{\ka}{\nu}\Bigr)\ad{\f}{\ka}\\
&\in\au{\omega}{\f} \Bigl(
\frac{A_{\omega\ka}}{\cont{\ka}{\omega}{\f}} \au{\ka}{\omega} +
\NN_{\omega\ka}^{\ka^-} \Bigr)\ad{\f}{\ka}
\end{split}\end{equation}


  Recall that the sets of operators
  $\{\geomie{\ka}{\theta}{\omega}{i}\}_i$ form a basis of
  $\NN_{\omega,\ka}$.
  A direct calculation shows that (compare \ref{c-g}, \ref{C-G})
 \begin{equation*} 
  \ad{\f}{\omega}\au{\ka}{\f}=
  \sum_{\theta\underset{i}{\hookrightarrow}\ka}
  \cont{\omega\oplus\coker(i)}{\f}{\F} \geomie{\ka}{\theta}{\omega}{i}
 \end{equation*}

  By the relations \ref{C-G} between the geometric and the
  cellular bases, we have that the coefficient of
  $\ceomie{\ka}{\ka}{\omega}{=}$ when expanding
  $\ad{\f}{\omega}\au{\ka}{\f}$ with respect to the cellular
  basis is given by
 \begin{align*}
&\sum_{\embi{\theta}{i}{\ka}} \cont{\omega \oplus\coker(i)}{\f}{\F}
\hat{\mu}_{\vec{\MM}/\vec{\TT}}
(\embi{\theta}{i}{\ka},\ka=\ka) \\
 =& \sum_{\embi{\theta}{i}{\ka}} \cont{\omega \oplus\coker(i)}{\f}{\F}
 \Bigr({\ka=\ka \atop \embi{\theta}{i}{\ka}}\Bigl)
 \chi(\coker(i))\\ 
 = &\sum_{\alpha}\sum_{i,~\coker(i)=\alpha}
 \cont{\omega \oplus \alpha}{\f}{\F}
  \Bigr({\ka=\ka \atop \embi{\theta}{i}{\ka}}\Bigl)
\chi(\alpha) \\ = &\sum_{\alpha} \cont{\omega \oplus
\alpha}{\f}{\F}
 \chi(\alpha)\sum_{i,~\coker(i)=\alpha}
 \Bigr({\ka=\ka \atop \embi{\theta}{i}{\ka}}\Bigl)\\
= &\sum_{\alpha} \cont{\omega \oplus \alpha}{\f}{\F}
\chi(\alpha)\Bigr({\ka \atop
 \alpha}\Bigl) && \text{(the duality axiom)} \\
= &\sum_{\alpha} \chi(\alpha) \hat{w}(\alpha) = s(0) =
 \avoid{\omega}{\f}{\ka}_{\F}&& \text{(equation~(\ref{specific}))}
\end{align*}

It follows that
\begin{equation}\label{last}\begin{split}
 \au{\omega}{\f}\Bigl(&\ad{\f}{\omega}\au{\ka}{\f}\Bigr)\ad{\f}{\ka} \\
&\in \au{\omega}{\f} \Bigl(
\avoid{\omega}{\f}{\ka}_{\F}\au{\ka}{\omega} +
\NN_{\omega\ka}^{\ka^-} \Bigr)\ad{\f}{\ka}
\end{split}\end{equation}

By Lemma~\ref{onto}, $\au{\omega}{\f}
\NN_{\omega\ka}^{\ka^-}\ad{\f}{\ka}<\HH_{\f}^{\ka^-}$, thus
comparing equations~(\ref{ddd}) and (\ref{last}) we get the
desired equation.

\end{proof}

\begin{Ex}[Fourier decomposition in the field case]

  Let $R=\Fq$ be the finite field with $q$ elements.
  Fix two natural numbers $m,n$ with $2m\le n$. Then $(n,m)$ is a symmetric couple. The group is $\GL_n(\Fq)$,
  the representation is $\FF_m$ and the Hecke
  algebra $\HH_m$. $(g_k)_{k\le m}$ and $(e_k)_{k\le m}$ are two
  bases for $\HH_m$.
  Example~\ref{final} gives for $k\le m$,
  \[ \gb_{k}=\sum_{i=k}^m\sum_{j=0}^i  (-1)^{i-k}q^{\bigl({i-k \atop
     2}\bigr)}\biggl({i \atop k}\biggr)_q
     \biggl( { n-i-j
     \atop m-i-j} \biggr)_{\! q}\cdot q^{j(m-i-j)}
     \biggl( {l-j \atop i-j} \biggr)_{\! q}
     \eb_j
  \]
(compare with \cite{dunkl}).
\end{Ex}

\section{Generalization of the theory} \label{general}

The theory developed in sections 2 and 3 is valid (after minor
changes of the terminology) for a large class of examples. In this
section we explain the necessary terminology, re-phrase some of the
theorems in a wider generality, and lastly, give a small list of
examples which fit into this framework.

\subsection{The general setting}

We begin by replacing the category of modules over a ring by an
arbitrary category. In order to speak about Grassmannians we need
the notion of a "subobject" in this generality. Our reference for
that is \cite[V\S7]{maclane}. For their fundamental importance in our
discussion we recall some of the definitions.

\begin{Def}
A morphism $i:x\rightarrow y$ in the category $\CC$ is called {\bf
monic} if for every object $z$ of $\CC$, the map
\[ i_*:\Hom(z,x)\rightarrow \Hom(z,y),\quad\phi\mapsto i\circ\phi \]
is injective.
\end{Def}

Let $\CC$ be a category. We abuse the notations and replace $\CC$
with its subcategory which has the same objects, but its morphisms
consist only the monics in $\CC$. We fix the category
$\CC$ for the rest of the section.

Let $\TT=\pi_0(\CC)$ denote the collection of types, that is
isomorphism classes, in $\CC$. We denote by
$\tau:\Ob(\CC)\rightarrow \TT$ the type map. Let
$*\!\longrightarrow\! *$ denote the category which consists of two
objects and one nonidentity arrow. Define $\vCC$ to be the
category of functors (and natural equivalences) from
$*\!\longrightarrow\! *$ to $\CC$. Observe that all the morphisms
in $\vCC$ are monics (by the assumption on $\CC$). We denote by
$\vTT$ the collection of types of $\vCC$ and by $\vec{\tau}$ the
type map.

\begin{Def}
Let $y$ be an object of $\CC$. A {\bf subobject} of $y$ is an
equivalence class of monics $x\overset{i}{\rightarrow} y$, under
the relation
\[ x\overset{i}{\rightarrow} y ~\sim~ x'\overset{i'}{\rightarrow} y \]
if and only if there exist an isomorphism
$x\overset{\phi}{\rightarrow} x'$ with $i'\circ\phi=i$.
\end{Def}

We denote by $\CC_y$ the collection of all subobjects of $y$. The
subobject of $y$ (represented by) $x\overset{i}{\rightarrow} y$ is
said to be smaller than the subobject (represented by)
$x'\overset{i'}{\rightarrow} y$ if there exist a (necessarily
monic) morphism $x\overset{\phi}{\rightarrow} x'$ with
$i'\circ\phi=i$. In our discussion we will assume that $\CC_y$ is
a {\bf finite lattice} for every $y\in\Ob(\CC)$. Symbols $\bigl({y
\atop \xi}\bigr)$ and $\bigl({y \atop \iota}\bigr)$ are, thus,
readily understood. These are $\Aut(y)$-spaces, hence yield to
representations of $\Aut(y)$ on the space of $\Q$-valued functions
defined on them.

\medskip

With the above terminology in hand, the reader is invited to observe
that sections 2 and 3 are valid {\em mutatis-mutandis} in this
generalized setting. In particular, we can define symmetric objects
and symmetric couples (see definitions~\ref{sym},\ref{symcoup}), and
deduce (in analogy with Theorem~\ref{multiplicities}):
\begin{The}
$F\in\Ob(\CC)$ be of symmetric type $\F$. For every $\la\le\F$,
denote by $\FF_{\la}$ the vector space of $\Q$-valued functions on
$\bigl({F \atop \la}\bigr)$. Let $\f$ be a type such that $\f \le
\F$ is a symmetric couple, and let $G=\Aut_{\CC}(F)$. There exists a
collection of non-equivalent irreducible $G$-representations
$\bigl\{\UU_{\la}\bigr\}_{\la \le \f}$ such that:
\begin{enumerate}
    \item $\FF_{\f}=\bigoplus_{\la \le \f}\UU_{\la}$.
    \item For every $\la,\nu \le \f$:
\[ \langle\UU_{\la},\FF_{\nu}\rangle=|\{\la \hookrightarrow \nu\}| \]
I.e., the multiplicity of $\UU_{\la}$ in $\FF_{\nu}$ is the number
of non-equivalent monics from an object of type $\la$ to an object
of type $\nu$. In particular $\UU_{\la}$ appears in $\FF_{\la}$
with multiplicity one and does not appear in $\FF_{\nu}$ unless
$\la \le \nu$.
\end{enumerate}

\end{The}

\subsection{Examples} \label{examples}

  We now give a list of examples of Categories.
  We will describe the symmetric objects
  and couples. We won't give proofs in all cases, as these can be
  regarded as easy exercises.
  We will try to describe the
  automorphism groups of symmetric objects, and their actions on the
  Grassmannians.

\subsubsection{Sets}

  In the category of finite sets, {\bf Sets}, the set of types
  $\TT_{{\bf Sets}}$ can be identified with $\N\cup\{0\}$, the
  type map given by $\tau(A)=|A|$. Every object in {\bf Sets} is
  symmetric, and the couple $(n,m)$ is symmetric \Iff\ $n\geq 2m$.
  For an object $A$ of type $n$, $\Aut(A)$ is identified with the
  permutation group $S_n$. The action of $\Aut(A)$ on $\bigl({A
  \atop m}\bigr)$ is identified with the action of $S_n$ on
  $\bigl({[n]\atop m}\bigr)$.

\subsubsection{Vector spaces} \label{Vec}

  Denote by ${\bf Vec}_q$ (or simply ${\bf Vec}$ when $q$ is given) the category
  of finite dimensional vector spaces over the finite field $\Fq$.
  $\TT_{{\bf Vec}}$ is naturally identified with $\N\cup\{0\}$,
  by letting for every $V\in \Ob({\bf Vec})$,
  $\ta(V)=\dim(V)$.
  Our notion of Grassmannian coincides with usual one.
  The cardinalities of the Grassmannians are given by the $q$-binomial
  functions:
  \[ \biggl({n \atop m}\biggr)_{\!{\bf Vec}} = \biggl({n \atop
     m}\biggr)_{\! q}=\frac{(q;q)_n}{(q;q)_m(q;q)_{n-m}},
  \]
  where $(a;q)_n$ is the Pochammer symbol and equals to
  $\prod_{i=0}^{n-1}(1-aq^i)$ .

  The automorphism groups $\Aut(V)$ is nothing but $\GL(V)$.
  It is easy to see that every type in $\TT_{{\bf Vec}}$ is
  symmetric, and that the couple $(n,m)$ is symmetric \Iff\ $n\geq 2m$ (this example is treated in \cite{dunkl}).

\subsubsection{Anti-symmetric bilinear forms} \label{Symp}

  Denote by ${\bf Symp}_q$ (or ${\bf Symp}$) the category which objects are couples
  $(V,B)$, where $V$ is an object of ${\bf Vec}_q$, and $B$ is a
  (possibly degenerate) anti-symmetric bilinear form on $V$ (we refer to such an object as a {\bf symplectic space}). The
  morphisms in ${\bf Symp}$ are given by
  \[ \Hom_{{\bf Symp}}((V,B),(V',B'))=\{ \phi\in\Hom_{{\bf
     Vec}}(V,V')~|~B=\phi^*(B')\}
  \]
  ${\bf Vec}$ appears as a full subcategory of ${\bf Symp}$, by
  $V\mapsto (V,0)$. The object of {\bf Vec} are
  denoted {\bf isotropic} when regarded as objects of {\bf Symp}.
  For the symplectic space $(V,B)$, we will use $i(V,B)$ to denote
  the dimension of a maximal isotropic subspace, that is
  \[ i(V,B)=\max\{\dim(U)~|~(U,0)<(V,B)\} \]
  We denote $\rad(V,B)=\{ v\in V~|~\forall u\in V,~B(u,v)=0 \}$,
  and $n(V,B)=\half\dim(V/\rad(V,B))$.
  Given $i=i(V,B)$ and $n=n(V,B)$, we have $\dim(V)=i+n$, and
  $\dim(\rad(V,B))=i-n$, thus the set of types may be identified
  with
  \[ \TT_{{\bf Symp}}=\{(i,n)~|~i,n\in\N\cup\{0\},~i\geq n\} \]
  The order structure on $\TT_{{\bf Symp}}$ is given by
  $(i,n)\geq (i',n')$ \Iff\ both $i\geq i'$ and $n\geq n'$.

  We already know that every isotropic object is symmetric.
  By Silvester theorem, if $B$ is non-degenerate, then
  $(V,B)$ is symmetric as well. On the other hand, it is easy to
  see that
  for $0<n<i$, the type $(i,n)$
  is non-symmetric (compare radical and non-radical lines).
  Therefore, we conclude that
  The symmetric types are exactly the types of the form
  $(i,0)$ or $(n,n)$. One easily verifies that the couples $((n,n),(n,0))$
are symmetric. The corresponding
  Grassmannians
  \[ \biggl( {(n,n) \atop (n,0)} \biggr) \]
  are known as the Lagrangian-Grassmannians.
  They are acted upon by $\Aut(n,n)=\Symp_n(\Fq)$.

\subsubsection{Bundles over sets} \label{Bun}

  Let $X$ be a set. We consider the category whose objects are the
  points of $X$, and which has a unique morphism between every two
  points in $X$. We denote this category by $X$ as well.
  Let $\mathcal{C}$ be a category.
  Let ${\bf Bun}=
  {\bf Bun}_X(\mathcal{C})$ be the category which objects are
  functors from $X$ to $\mathcal{C}$, and a morphism between to
  such functors, $E$ and $F$, is a couple $(\sigma,\phi)$, where
  $\sigma$ is a permutation of $X$ and $\phi$ is a natural
  transformation from $E$ to $F\circ\sigma$.
  The types of ${\bf Bun}$ are given by the
  multiplicities of the various types in $\mathcal{C}$, that is
  \[ \TT_{{\bf Bun}}=\{ r:\TT_{\mathcal{C}}\rightarrow
     \N\cup\{0\}~|~\sum_{\lambda\in\TT_{\mathcal{C}}}r(\lambda)=|X|
     \}
  \]
  in particular we denote by $\tilde{\lambda}$ the type in
  $\TT_{{\bf Bun}}$ which satisfies
  $\tilde{\lambda}(\lambda)=|X|$. It is immediately seen that an object $E$ in ${\bf Bun}$ is
  symmetric \Iff\ $E_x$ is symmetric for every $x\in X$.
  Equivalently, the symmetric types of ${\bf Bun}$ are those
  supported on the set of symmetric types of $\mathcal{C}$.
  For a symmetric couple in $\mathcal{C}$, $(\F,\f)$,
  $(\tilde{\F},\tilde{\f})$ is a symmetric couple in ${\bf Bun}$. The group $\Aut(\tilde{\F})$ can be identified with the wreath
  product $S_{|X|}\propto\Aut(\F)$. It acts on the Grassmannian
  $\bigl({\tilde{\F} \atop \tilde{\f}}\bigr)$, which, in turn,
  can be identified with $\bigl({\F \atop \f}\bigr)^{|X|}$.  Important special cases occur when $\mathcal{C}$ is
  {\bf Sets}, {\bf Vec} or {\bf Symp}. For short we will consider
  in the sequel only $\mathcal{C}={\bf Vec}$. The other two cases
  are similar.

\subsubsection{Vector bundles} \label{VBun}

  Let ${\bf VBun}_q={\bf Bun}_X({\bf Vec}_q)$ be the category of finite
  dimensional $\Fq$-vector bundles over the set $X$. The types in ${\bf VBun}$ are given by partitions of $k=|X|$,
  $r=(r_1,\ldots,r_m)$, where
  \[ r_i=|\{ x\in X~:~ \dim(E_x)=i\}| \]
  The ordering is given by
  \[  r \geq r'~~\Leftrightarrow~~
      \forall j,~\sum_{i\geq j}r_i\geq \sum_{i\geq
      j}r'_i
  \]
  Another identification of the set of types is
  \[ \TT_{{\bf VBun}}=\{ \la=(\la_1,\la_2,\ldots,\la_k)~|~\la_1
     \geq \la_2\geq\cdots\geq\la_k\ge 0 \}
  \]
  where the type map is given by
  \[ \tau(E)=\la(E)=(\dim(E_{x_1}),\dim(E_{x_2}),\ldots,
     \dim(E_{x_k}))
  \]
  Here we assume $X=\{x_1,\ldots,x_k\}$ is an ordering of $X$,
  such that $\dim(E_{x_i})\geq\dim(E_{x_{i+1}})$.
  The ordering is now given by
  the lexicographic order
  \[ \la\geq\la'~~\Leftrightarrow~~\forall i,~\la_i\geq\la'_i \]
  Using these coordinates, we rewrite $\tilde{m}=(m,m,\ldots,m)$.
  Denote
  \[ \TT_{\leq\tilde{m}}=\{\la\in\TT_{{\bf
     VBun}}~|~\la\le\tilde{m} \}
  \]
  \[ \Lambda^m_k=\{
     (\la_1,\ldots,\la_k)~|~m\ge\la_1\ge\cdots\ge\la_k\ge 0 \}
  \]
  Then $\TT_{\le\tilde{m}}=\Lambda^m_k$. Let $m,n$ be given, and assume $n\ge 2m$. The couple
  $(\tilde{n},\tilde{m})$ is a symmetric couple.
  The automorphism group of an object of type $\tilde{n}$ is
  isomorphic to the wreath product $S_k \propto \GL_n(\Fq)$. It
  acts on the corresponding Grassmannian of $\bigl({n \atop
  m}\bigr)^k$ elements.



\vspace{\bigskipamount} \vspace{\bigskipamount}
\vspace{\bigskipamount}

\begin{footnotesize}
\begin{quote}

Uri Bader\\
Department of Mathematics, University of Chicago,\\
Chicago, Illinois 60637, USA \\
{\tt bader@math.uchicago.edu} \\ \\

Uri Onn\\
Th\'eorie des Groupes,
U.F.R.\ de Math\'ematiques, Universit\'e Paris 7,\\
Case 7012,$\;$ 2, place Jussieu, 75251 Paris cedex 05, France \\
{\tt onn@math.jussieu.fr} \\ \\

\end{quote}
\end{footnotesize}

\end{document}